# ON THE VARIATIONAL DISTANCE OF TWO TREES[1]

BY M. A. STEEL AND L. A. SZÉKELY[2]

*University of Canterbury and University of South Carolina*

A widely studied model for generating sequences is to "evolve" them on a tree according to a symmetric Markov process. We prove that model trees tend to be maximally "far apart" in terms of variational distance.

**1. Introduction.** In this paper we investigate sequences that have been generated on the tree by a simple Markov model. Such processes are widely-studied in molecular genetics, and in other areas of applied probability (including broadcasting and statistical physics). More precisely, we study the separation—as measured by variational distance—of the probability distribution on sequence patterns generated by different trees. We find that a large tree generates a probability distribution that is typically at maximal distance from that generated by nearly all other trees.

To describe our results more precisely, we first provide some terminology concerning trees and random processes on them. In a tree, vertices of degree 1 are called *leaves*, as opposed to *internal vertices*. A tree is *binary* if all vertices have degree 1 or 3. Consider a set $X$ of labels. A *phylogenetic $X$-tree* is a tree in which leaves are identified with elements of $X$. (We do not require phylogenetic $X$-trees to be binary by definition for technical reasons, as we will have to deal with subtrees of phylogenetic $X$-trees.) We will regard two phylogenetic $X$-trees as being identical if there is a graph isomorphism between them, which, in addition, if restricted to $X$, is the identity function of $X$. If $|X| = n$, then the number of different binary phylogenetic $X$-trees is $(2n-5)!! \ [= 1 \times 3 \times 5 \times \cdots \times (2n-5)]$ [16]. For a phylogenetic $X$-tree

Received August 2005; revised January 2006.
[1]Supported in part by the NZIMA (Maclaurin Fellowship).
[2]Supported in part by NSF Grants DMS-00-72187 and DMS-03-02307.

$\mathcal{T}$, let $[\mathcal{T}]$ denote the corresponding unlabeled tree. The *distance* $d_\mathcal{T}(u,v)$ between two vertices, $u,v$ in a tree $\mathcal{T}$ is the number of edges on the unique path connecting them.

We now describe a model for the evolution of *binary* sequences on a tree. This model has been described by various authors (and in a range of disciplines, including molecular biology, information theory and physics; for references, see [8, 16]). Here we refer to this model as the *CFN model* (short for Cavender–Farris–Neyman model); it has also been referred to in the literature as the "symmetric binary channel" and the "symmetric 2-state Poisson model." The CFN model provides a simple model for the evolution of purine–pyrimidine sequences. The significance of this simple model is that phenomena shown for the CFN model often extends to more realistic models of sequence evolution, and we will describe how our main results concerning the CFN model generalize. The term *CFN tree* will refer to a phylogenetic $X$-tree equipped with a CFN model.

Suppose we have two states, 0 and 1, and a phylogenetic $X$-tree $\mathcal{T}$. The CFN model assigns probabilities to the patterns of state of the elements of $X$ as follows. Let us associate a number $p_e$ ($0 < p_e < 1/2$) with the edge $e$ called the *transition probability*. Let $\xi_e$ denote a random indicator variable associated to edge $e$ with $\mathbb{P}[\xi_e = 1] = p_e$, and assume the $\xi_e$'s are independent. Fix any vertex $v$ and assign state 0 or 1 to $v$ with equal probability $1/2$. Note that, for every vertex $u$ of $\mathcal{T}$, there is a unique path denoted $path(u,v)$ in $\mathcal{T}$ and so we may define

$$(1.1) \qquad state(u) = state(v) + \sum_{e \in path(u,v)} \xi_e \bmod 2.$$

This gives a (joint) probability distribution on the set of all assignment of states (0 or 1) to the vertices of $\mathcal{T}$, and thereby a marginal distribution on state assignments to the leaves of $\mathcal{T}$—we call each such assignment $\chi: X \to \{0,1\}$ a (state) *pattern*, and we let $\mathcal{P}_\chi$ denote the probability of generating $\chi$ under this model.

The CFN model is thus specified by the pair $(\mathcal{T}, \mathcal{P})$, where $\mathcal{P}$ is the map that associated to each edge $e$ its transition probability $p(e)$. We refer to $\mathcal{T}$ as a *CFN tree* and $\mathcal{P}$ as a *transition mechanism*.

The probability $p$ that the endpoints of a path $uw$ in a CFN tree $\mathcal{T}$ are in different states is nicely related to the transition probabilities of edges of the $uw$-path:

$$(1.2) \qquad p = \tfrac{1}{2}\left(1 - \prod_{e \in path(u,w)}(1-2p_e)\right).$$

Formula (1.2) is well known and is easy to prove by induction. Formula (1.2) also shows that the transition probability of a path is not less than the largest transition probability on its edges. It is well known [18] that (1) changing the location of $v$ in $\mathcal{T}$, or (2) substituting a path by a *single edge* in a CFN tree, and assigning to the new edge a transition probability according to (1.2) *does not change* the probability distribution of patterns.

Usually $k$ independent experiments are made to generate random patterns from a binary CFN tree $\mathcal{T}$, they are called *sites*. The (abstract) *phylogeny reconstruction problem* is the following: from the observed pattern frequencies, determine, with a prescribed probability, what was the underlying binary phylogenetic $X$-tree. We have shown in [6] that if $|X| = n$ and $n \to \infty$, then $k = \Omega(\log n)$ sites are needed to return the true underlying tree with probability at least $\frac{1}{2} + \varepsilon$ with either a deterministic algorithm or with a randomized algorithm whose random bits are independent from the random events on the CFN tree. Sequence length requirements for accurate tree reconstruction is not only of mathematical interest, but also a topical issue in molecular systematics (e.g., [3, 15]). We showed in [6] that, for fixed $0 < f \leq g < 1/2$, $f \leq p_e \leq g$, and $n \to \infty$, phylogeny reconstruction is possible for all model trees, when $k$ is a certain polynomial of $n$; is possible for some model trees, when $k$ is a logarithmic function of $n$; and is possible for almost all model trees, either in the uniform random binary $X$-tree model or in the Yule–Harding model, when $k$ is a certain polylogarithmic function of $n$. More recent work by Mossel and colleagues [5, 12] has established further instances for which logarithmic dependence of $k$ on $n$ suffices for accurate tree reconstruction and cases for which polynomial dependence is necessary.

In this paper we show *asymptotic results*. The theorems are about $n$-leaf trees, but their conclusions are $o(1)$ (limit) relations as $n \to \infty$. The understanding is that, for a *sequence* of $n$-leaf trees satisfying the hypotheses, the limit relation holds. It would be technically more proper to speak about *sequences of trees* in the statements of the theorems, but we follow the tradition of random graph theory [1, 4] not speaking explicitly about sequences. With the exception of Section 4, we study problems where the bounds on $p_e$ are *fixed*, and we let $n \to \infty$. In Section 4 we show that many of the results generalize if dependence of the bounds on $n$ is allowed but limited.

**2. Results.** Let us be given two binary phylogenetic $X$-trees $\mathcal{T}_1, \mathcal{T}_2$ with CFN transition mechanism $\mathcal{P}_1$ and $\mathcal{P}_2$, respectively. The variational distance of their pattern distributions is

$$(2.1) \qquad \text{vardist}((\mathcal{T}_1, \mathcal{P}_1),\ (\mathcal{T}_2, \mathcal{P}_2)) = \sum_{\chi} |(\mathcal{P}_1)_\chi - (\mathcal{P}_2)_\chi|.$$

This distance lies between 0 and 2, and in Theorem 3.1 we show that *almost all* binary trees are maximally distant (in terms of variational distance) from any given binary tree with a given CFN transition mechanism, under mild assumptions on their transition mechanisms. A practitioner may argue that Theorem 3.1 has limited relevance, since the uniform distribution of trees is just one particular prior distribution on trees, and the CFN model is very particular. However, the conclusion of Theorem 3.1 holds not just for the counting measure, but for all permutation invariant measures on phylogenetic $X$-trees; moreover, it holds for more general, and for the applications more realistic classes of transition mechanisms (Theorem 4.1). This result may not be surprising: as we equip randomly selected trees with CFN models, they have many local statistics that are essentially independent and have different marginals in the two trees. Therefore, analogously to the Kakutani dichotomy, their measures are expected to be (near) orthogonal.

Farach and Kannan [2, 9, 10] designed an algorithm for phylogeny reconstruction based on convergence to the true tree in variational distance and suggested to pay more attention to the variational distance in phylogeny reconstruction. Some support for the utility of this metric is provided by results that we present in Sections 3 and 4: if we get just *close* to a model tree in variational distance, then we already excluded most of the false candidates for the phylogenetic tree.

However, a simple fact provides a sharp contrast to the results mentioned above. Note that in practice we estimate the model distribution of patterns by the *observed frequency* of patterns. For *sub-exponential* sequence length, which is known to be sufficient for phylogeny reconstruction with probability $1 - o(1)$ as $0 < f \leq g < 1/2$ fixed and $f \leq p_e \leq g$, as $n \to \infty$ (see the discussion in Section 1), the variational distance between the model pattern distribution and the observed pattern distribution is near 2 with probability $1 - o(1)$. (For details, see our technical report [19].)

In other words, phylogeny reconstruction is well possible *without* convergence of the observed pattern distribution to the model pattern distribution in variational distance.

Therefore, the accuracy of tree reconstruction cannot be captured by variational distance alone. This conclusion was suggested by [7] and [14], though with less explicit theoretical justification.

## 3. Variational distance of CFN trees is typically large.

THEOREM 3.1. *Fix $0 < f$ and $g < 1/2$. There exists a function $\varepsilon(n) = \varepsilon_{f,g}(n) = o(1)$ as $n \to \infty$, such that, for every binary phylogenetic $X$-tree $\mathcal{T}_1$ with CFN transition mechanism $\mathcal{P}_1$ where $p_e \leq g$ in $\mathcal{P}_1$, the following holds: For almost all [i.e., $(1 - o(1))(2n - 5)!!$ in number] binary phylogenetic $X$-trees $\mathcal{T}_2$, equipped with an arbitrary transition mechanism $\mathcal{P}_2$, where $f \leq p_e$ in $\mathcal{P}_2$, we have*

$$\text{vardist}((\mathcal{T}_1, \mathcal{P}_1), (\mathcal{T}_2, \mathcal{P}_2)) \geq 2 - \varepsilon(n). \tag{3.1}$$

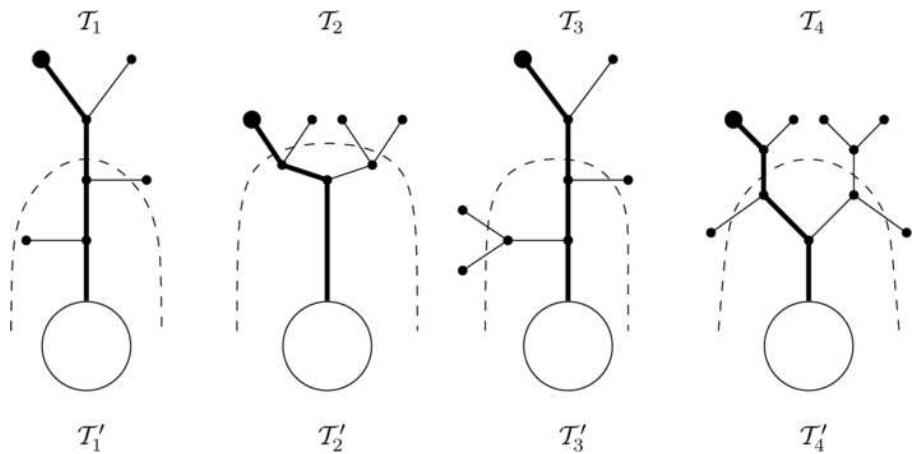

Fig. 1. *Ending of a longest path in a binary tree.*

The proof requires a number of lemmas, which we now state.

LEMMA 3.2. *For every binary phylogenetic $X$-tree $\mathcal{T}$ on $n \geq 4$ leaves, there are at least $n/4$ disjoint pairs of leaves $a_i, b_i$, such that, for every $i$:*

(i) *$a_i$ and $b_i$ are separated by a distance of 2 or 3;*
(ii) *for $i \neq j$, the $a_i b_i$ and the $a_j b_j$ paths in $\mathcal{T}$ are edge disjoint.*

PROOF. The claim is true for $4 \leq n \leq 8$, since then any longest path ends in two disjoint cherries. This is the basis for an induction proof on $n$. It is easy to see that, for $n \geq 9$, there exists a longest path in $\mathcal{T}$, for which one end must be a leaf in a cherry that lies at the top portion of the tree given by one of the four cases shown in Figure 1 (the other end of the path lies in the bottom part of the tree, represented by a circle). In each of the four cases truncate the tree $T_i$ as indicated by the dashed curve to obtain $T'_i$. For $i = 1, 2, 3, 4$, $T'_i$ has $n-2$ (resp. $n-2$, $n-3$, $n-4$) leaves, and the induction hypothesis applies to $T'_i$. In all four cases it is easy to add two new close vertex pairs to create the required set of them for $T_i$, while destroying at most one which pre-existed in $T'_i$. □

REMARK 3.3. As Figure 2 shows, the conclusion of Lemma 3.2 is essentially the best possible.

LEMMA 3.4 (Tree-chopping lemma, [17], Lemma 3). *Let $\mathcal{T}$ be an arbitrary binary $X$-tree and $q \geq 2$ integer. Then edges can be deleted from $\mathcal{T}$ such that a forest results with the following properties:*

(i) *The number of leaves from $X$ in any tree of the forest is at most $2q - 2$.*

(ii) *The number of leaves from $X$ in any tree of the forest is at least $q$, except possibly for one tree. (We shall call this exceptional tree* degenerate.*)*

Recall the Azuma–Hoeffding inequality (see [1]):

LEMMA 3.5. *Suppose $\mathbf{X} = (X_1, X_2, \ldots, X_k)$ are independent random variables taking values in any set $S$, and $L:S^k \to \mathbb{R}$ is any function that satisfies the condition: $|L(\mathbf{u}) - L(\mathbf{v})| \leq t$ whenever $\mathbf{u}$ and $\mathbf{v}$ differ at just one coordinate. Then,*

$$\mathbb{P}[|L(\mathbf{X}) - \mathbb{E}[L(\mathbf{X})]| \geq \lambda] \leq 2\exp\left(-\frac{\lambda^2}{2t^2 k}\right). \tag{3.2}$$

The following lemma is obvious.

LEMMA 3.6. *Let $\mathcal{F}$ denote a fixed phylogenetic $X$-tree, with $|X| = n$, and let $\tau = [\mathcal{F}]$ (the corresponding unlabeled tree). Let $\pi$ be a randomly selected permutation of $X$ under the uniform distribution. Let $\pi(\mathcal{F})$ denote the phylogenetic $X$-tree that we obtain from $\mathcal{F}$ by changing all leaf labels from $v$ to $\pi(v)$ simultaneously. Then $\pi(\mathcal{F})$ represents a random uniform selection from those binary phylogenetic $X$-trees whose underlying unlabeled tree is $\tau$.*

From now on, for notational convenience, we pretend that 4 divides $n$.

LEMMA 3.7. *For an $X$ with $|X| = n$, and $n/4$ disjoint $a_i, b_i$ ordered pairs from $X$, there exist functions $m(n) \to \infty$, $h(n) \to \infty$ and $g(n) \to \infty$, such that the following holds. For every unlabeled binary tree $\tau$ with $n$ leaves, for all but a $\frac{1}{g(n)}$ fraction of binary phylogenetic $X$-trees $\mathcal{T}$ with property $[\mathcal{T}] = \tau$, there is an index set $I$ such that $|I| = m(n)$ and:*

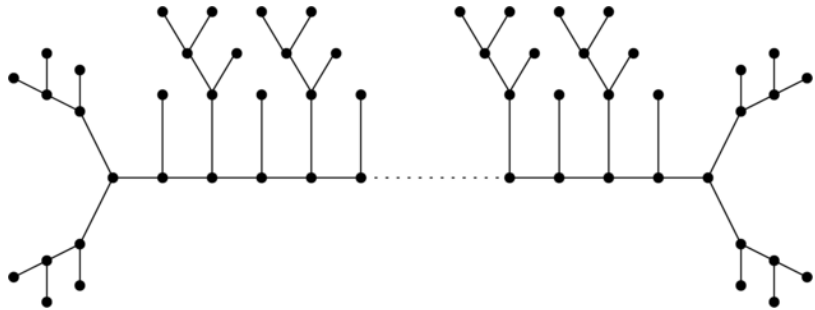

FIG. 2. *Binary tree on $4t + 9$ leaves, with only $t + 3$ close leaf pairs.*

  (i) $d_\mathcal{T}(a_i, b_i) \geq h(n)$ *for all* $i \in I$; *and*
  (ii) *for* $i, j \in I$, $i \neq j$, $path_\mathcal{T}(a_i, b_i)$ *and* $path_\mathcal{T}(a_j, b_j)$ *are edge disjoint.*

PROOF. Let $\mathcal{F}$ denote a fixed binary phylogenetic $X$-tree such that $[\mathcal{F}] = \tau$, with $|X| = n$. Apply Lemma 3.4 to $\mathcal{F}$ with $q = \lceil \log^2 n \rceil$. Let $L_1, L_2, \ldots, L_s$ denote the leaf sets that the nondegenerate trees contain from $X$. From the lemma, $q \leq |L_i| \leq 2q - 2$, and at most $q - 1$ elements of $X$ are not in some $L_i$. Let $\pi$ be a randomly selected permutation of $X$ under the uniform distribution. Let $\pi(\mathcal{F})$ denote the phylogenetic $X$-tree that we obtain from $\mathcal{F}$ by changing all leaf labels from $v$ to $\pi(v)$ simultaneously. According to Lemma 3.6, $\pi(\mathcal{F})$ represents a random uniform selection from those binary phylogenetic $X$-trees whose underlying unlabeled tree is $\tau$. The previous application of Lemma 3.4 still partitions $\pi(\mathcal{F})$, the leaf sets of the nondegenerate trees intersect $X$ in $\pi(L_1), \pi(L_2), \ldots, \pi(L_s)$, and we still have $q \leq |\pi(L_i)| \leq 2q - 2$. Therefore, for $i \neq j$, if $path^i_{\pi(\mathcal{F})}$ (resp., $path^j_{\pi(\mathcal{F})}$) connects an arbitrary pair of vertices of $L_i$ (resp., $L_j$) in the tree $\pi(\mathcal{F})$, then

$$(3.3) \qquad path\ i^i_{\pi(\mathcal{F})} \text{ is edge disjoint from } path^j_{\pi(\mathcal{F})}.$$

Set $h(n) = \log \log n$ and $m(n) = \frac{n}{4q-2} - \frac{1}{2}$. Observe from Lemma 3.4 and the choice of $q$ that $n \leq (s+1)(2q-2)$, and therefore,

$$(3.4) \qquad m(n) \leq \frac{s}{2}.$$

We are going to find an appropriate $g(n)$ for this choice. We call a leaf set $Y \subset X$ *infected*, if there is a $1 \leq j \leq n/4$, such that *both* $a_j, b_j \in Y$. Let $E$ denote the event that, for our fixed $\tau$ and $\mathcal{F}$, $\pi(\mathcal{F})$ has the property that for all $j = 1, 2, \ldots, s$, $\pi(L_j)$ is infected; and let $F$ denote the event that, in addition to $E$, for at least half of the indices $j = 1, 2, \ldots, s$, one finds some $i_j$, such that both $a_{i_j}, b_{i_j} \in \pi(L_j)$ [i.e., they do infect $\pi(L_j)$] and $d_{\pi(\mathcal{F})}(a_{i_j}, b_{i_j}) \geq h(n)$. In view of (3.3), the $a_{i_j}, b_{i_j}$ paths in $\pi(\mathcal{F})$ are pairwise edge disjoint for $j = 1, 2, \ldots, s$.

Observe that

$$(3.5) \qquad \mathbb{P}[\pi(L_j) \text{ not infected}] = \frac{\sum_{u=0}^{|L_j|} \binom{n/4}{u} 2^u \binom{n/2}{|L_j|-u}}{\binom{n}{|L_j|}}.$$

[A noninfected $L_j$ can have zero or one element from every $(a_i, b_i)$ pair, for $i = 1, 2, \ldots, n/4$. The case analysis is based on the number $u = |\pi(L_j) \cap \{a_i, b_i : i = 1, 2, \ldots, n/4\}|$. There are $\binom{n/4}{u}$ to select a subset of $u$ indices from $\{1, 2, \ldots, n/4\}$, and then $2^u$ ways to tell if $a_i$ or $b_i$ selected for the particular index set into $L_j$. There are $\binom{n/2}{|L_j|-u}$ ways to make $L_j$ complete using $|L_j| - u$ elements not belonging to $\{a_i, b_i : i = 1, 2, \ldots, n/4\}$.]

Comparison of consecutive terms show that the largest term in the numerator of the RHS of (3.5) is $u = |L_j|$. Using the usual notation $(x)_m$ for the $m$th falling factorial, it follows that

$$\text{(3.6)} \quad \mathbb{P}[\pi(L_j) \text{ not infected}] \leq \frac{(|L_j|+1)2^{|L_j|}\binom{n/4}{|L_j|}}{\binom{n}{|L_j|}}$$

$$\text{(3.7)} \quad = \frac{(n/4)_{|L_j|}}{(n)_{|L_j|}}(|L_j|+1)2^{|L_j|}$$

$$\leq \frac{n^{|L_j|}}{4^{|L_j|}(n-|L_j|)^{|L_j|}}(|L_j|+1)2^{|L_j|}$$

$$\text{(3.8)} \quad \leq (1+o(1))2^{-|L_j|}(|L_j|+1)$$

$$\leq (1+o(1))2^{-q}(2q-1) \leq 2^{-0.01\log^2 n},$$

and from (3.6)–(3.8),

$$\text{(3.9)} \quad \mathbb{P}[\exists j : \pi(L_j) \text{ not infected}] \leq \frac{n}{q} 2^{-0.01\log^2 n}.$$

By (3.9), we showed that

$$\text{(3.10)} \quad \mathbb{P}[E] > 1 - n2^{-0.01\log^2 n}.$$

Call the ordered $s$-tuple of pairwise disjoint sets $Y_1, Y_2, \ldots, Y_s \subset X$ *feasible*, if $|Y_i| = |L_i|$ and $Y_i$ is infected for $i = 1, 2, \ldots, s$. Now we turn to the conditional probability $\mathbb{P}[F|E]$. Observe

$$\text{(3.11)} \quad \mathbb{P}[F|E] = \sum_{Y_1,Y_2,\ldots,Y_s \text{ feasible}} \mathbb{P}[F|\forall i : \pi(L_i) = Y_i]\mathbb{P}[\forall i : \pi(L_i) = Y_i]$$

$$\text{(3.12)} \quad \leq \max_{Y_1,Y_2,\ldots,Y_s \text{ feasible}} \mathbb{P}[F|\forall i : \pi(L_i) = Y_i].$$

Assume now that an arbitrary feasible $Y_1, Y_2, \ldots, Y_s$ is *fixed*. A $\pi$ that satisfies the condition in (3.12) is nothing else but the juxtaposition of $\pi_i :_i \to Y_i$ bijections for $i = 1, 2, \ldots, s+1$. Therefore, a uniform random $\pi$ satisfying the condition in (3.12) can be realized by a sequence of *independent* uniform random choices of bijections $\pi_i$ from $L_i$ to $Y_i$, $i = 1, 2, \ldots, s+1$.

Let $\pi_i : L_i \to Y_i$ denote a uniform random bijection for $i = 1, 2, \ldots, s+1$. Conditional on $E$, for every $i = 1, 2, \ldots, s$, fix an $a_{i_j}, b_{i_j}$ leaf pair that infects $Y_i$. Observe that the conditional event

$$F|\forall i : \pi(L_i) = Y_i$$

is implied, if for at least half of the indices $1 \leq i \leq s$, we have $d_{\pi(\mathcal{F})}(a_{i_j}, b_{i_j}) \geq h(n)$. Also observe, that notwithstanding the notation $d_{\pi(\mathcal{F})}$, this distance depends *only* on the *single* $\pi_i$ under consideration. No matter what is the value of $\pi_i^{-1}(a_{i_j})$, at most $2^{h(n)}$ vertices of $L_i$ can be closer than $h(n)$ to $\pi_i^{-1}(a_{i_j})$ in the binary tree $\mathcal{F}$. Those at most $2^{h(n)}$ vertices can be pre-images of $b_{i_j}$ under $\pi_i$ (and $\pi$ as well), if $d_{\pi(\mathcal{F})}(a_{i_j}, b_{i_j}) < h(n)$. Therefore,

$$\mathbb{P}[d_{\pi(\mathcal{F})}(a_{i_j}, b_{i_j}) \geq h(n)] \geq 1 - \frac{2^{h(n)}}{|L_i|} = 1 - \frac{2^{\log\log n}}{\log^2 n} = 1 - \frac{1}{\log^{2-\log 2} n}.$$

Hence, a lower bound for $\mathbb{P}[F|E]$ is the probability of at least $s/2$ successes in a sequence of $s$ independent Bernoulli trials, each with probability of success $p = 1 - \frac{1}{\log^{2-\log 2} n}$. Not having at least $m(n)$ successes implies not having at least $s/2$ successes by (3.4), and probability of the latter event can easily be bounded from above by Lemma 3.5 ($t=1$, $k=s$, $\lambda=s/3$), as soon as $\frac{1}{\log^{2-\log 2} n} < 1/6$, by

(3.13) $$2e^{-s/18}.$$

Finally, using (3.10) and (3.13), we have

(3.14) $$1 - \mathbb{P}[F] = 1 - \mathbb{P}[E] + \mathbb{P}[E](1 - \mathbb{P}[F|E])$$
$$\leq n 2^{-0.01\log^2 n} + 2e^{-n/(64\log^2 n)},$$

and since the RHS of (3.14) is $o(n)$, we can take for $g(n)$ its reciprocal. $\square$

PROOF OF THEOREM 3.1. Specify now $n/4$ leaf pairs $\{a_i, b_i\}$ of $\mathcal{T}_1$ according to Lemma 3.2—for notational convenience, we assume again that $n$ is a multiple of 4. We set $m(n)$, $h(n)$, $g(n)$ and $I$ according to the statement of Lemma 3.7. We are going to show that, for every fixed $(\mathcal{T}_1, \mathcal{P}_1)$ and fixed unlabeled tree $\tau$, if $[\mathcal{T}_2] = \tau$ and $\mathcal{T}_2$ is not in the exceptional set of trees described in Lemma 3.7, then the variational distance between $(\mathcal{T}_1, \mathcal{P}_1)$ and $(\mathcal{T}_2, \mathcal{P}_2)$ differs from 2 by at most a quantity that is $o(1)$ as a function of $n$. Recall that $state(x)$ denotes the state of leaf $x \in X$ in a CFN tree. Consider the random indicator variable $Z_i$, which is 1, if $state(a_i) = state(b_i)$, and 0 otherwise, and $Z = \sum_{i \in I} Z_i$, which depends on the distribution of leaf colorations of the CFN tree. We will speak about $Z_i^{(1)}$, $Z^{(1)}$ and $Z_i^{(2)}$, $Z^{(2)}$ as the CFN tree is $(\mathcal{T}_1, \mathcal{P}_1)$ or $(\mathcal{T}_2, \mathcal{P}_2)$, and similarly about $state_1$ and $state_2$, and will drop the superscript if the argument applies to both.

By the linearity of expectation,

(3.15) $$\mathbb{E}[Z] = \sum_{i \in I} \mathbb{E}[Z_i] = \sum_{i \in I} \mathbb{P}[state(a_i) = state(b_i)].$$

In $(\mathcal{T}_1, \mathcal{P}_1)$, we have $\mathbb{P}[state_1(a_i) \neq state_1(b_i)] \leq \frac{1}{2}(1 - (1-2g)^3)$, by (1.2), and hence,

(3.16) $$1 - 3g + 6g^2 - 4g^3 \leq \mathbb{P}[state_1(a_i) = state_1(b_i)].$$

Formula (3.15) and inequality (3.16) imply that

$$\mathbb{E}[Z^{(1)}] \geq (1 - 3g + 6g^2 - 4g^3)m(n) \tag{3.17}$$

In $(\mathcal{T}_2, \mathcal{P}_2)$, by a similar argument, we have

$$\mathbb{P}[state_2(a_i) = state_2(b_i)] \leq 1 - \tfrac{1}{2}(1 - (1-2f)^{h(n)}) = \tfrac{1}{2} + o(1) \tag{3.18}$$

by (1.2), and $h(n) \to \infty$. By linearity (3.15), we have

$$\mathbb{E}[Z^{(2)}] \leq (1 + o(1))\frac{m(n)}{2}. \tag{3.19}$$

We are going to show that, with high probability, both $Z^{(1)}$ and $Z^{(2)}$ are very close to their respective expectations. This will be easy to show, since both of them are the sums of independent indicator variables. [Use Lemma 3.5 for $X_i = Z_i^{(1)}$ (resp. $Z_i^{(2)}$), $k = m(n)$, $t = 1$, $\lambda = m(n)^{2/3}$.]

It is easy to see that, for $0 < g < 1/2$, we have

$$1/2 < 1 - 3g + 6g^2 - 4g^3, \tag{3.20}$$

and therefore, using (3.17) and (3.19), $\mathbb{E}[Z^{(1)}]$ and $\mathbb{E}[Z^{(2)}]$ are separated by a linear function of $m(n)$, for example, $l(n) = \tfrac{1}{2}(1 - 3g + 6g^2 - 4g^3 + \tfrac{1}{2})m(n)$. Consider now the event $H$: "$Z > l(n)$." In $(\mathcal{T}_1, \mathcal{P}_1)$, event $H$ has probability $1 - o(1)$, while in $(\mathcal{T}_2, \mathcal{P}_2)$, the complement of event $H$ has probability $1 - o(1)$. This implies that the variational distance of $(\mathcal{T}_1, \mathcal{P}_1)$ and $(\mathcal{T}_2, \mathcal{P}_2)$ is $2 - o(1)$. □

**4. Variational distance in more general models.** In this section we provide a result (Theorem 4.1) that is a three-fold generalization of Theorem 3.1. The three extensions allow (i) more general probability distributions on trees ("permutation-invariant measures"), (ii) more general transition models than the CFN model ("conservative, separable processes") and (iii) a weakening of the constraints on the parameters of the model.

*Permutation-invariant measures on trees.* Let us call a measure $\mu$ on the set of $(2n-5)!!$ binary phylogenetic $X$-trees *permutation invariant*, if for every $\pi$ permutation of $X$ and any phylogenetic $X$-tree $\mathcal{F}$, $\mu(\mathcal{F}) = \mu(\pi(\mathcal{F}))$. Note that Lemma 3.6 stated that the uniform distribution (or counting measure) on binary phylogenetic $X$-trees is permutation invariant. A practitioner may argue that Theorem 3.1 has limited relevance, since the uniform distribution of trees is just one particular prior distribution on trees. However, any relevant distribution of trees is permutation invariant and it is easy to see that the stronger Theorem 3.1 holds with basically the same proof. A nonuniform, phylogenetically relevant permutation invariant distribution on phylogenetic $X$-trees is the *unrooted Yule–Harding distribution* [6].

*More general transition processes* (*conservative, separable processes*). The restriction of the CFN to two states and symmetric transition probabilities is convenient for description and proofs. However, much of the argument used in the proof of Theorem 3.1 can be generalized to models that are much closer to those used in modern molecular biology. We identify two key properties that are used in the proof, and that both apply to a range of substitution models.

Suppose we have a set $S$ of $q \geq 2$ states. A *pattern* will now refer to a state assignment function $\chi : X \longrightarrow S$, where $X$ is the leaf set of $\mathcal{T}$. Assume that we have a probability distribution on the patterns of a binary phylogenetic $X$-tree, where $\mathcal{P}_\chi$ denotes the probability of pattern $\chi$. Selecting a random pattern according to the distribution, we can observe a random *state* of any particular leave. For a pair of leaves $a, b$, let $E(a, b)$ be the event that $state(a) = state(b)$. Let us be given a strictly decreasing function $H : [0, \infty) \to (c, 1]$ with $H(0) = 1$, and a $c > 0$ constant, such that $\lim_{x \to \infty} H(x) = c$. We assume that $H$ and $c$ are fixed and do not depend on $n$. We say that a probability distribution on patterns is *conservative* if

(C)   there exists an assignment of $t(e) > 0$ to each edge $e$ of $\mathcal{T}$,
so that the following condition holds: For each pair $a, b \in X$,
$\mathbb{P}[E(a,b)] = H(\sum_{e \in path(a,b)} t(e))$.

The CFN model satisfies condition (C), as can easily be seen from (1.2) by taking $t(e) = -\frac{1}{2}\log(1 - 2p_e)$, $H(x) = \frac{1}{2}(1 + \exp(-2x))$, and $c = \frac{1}{2}$. More generally, condition (C) is satisfied by any tree-based Markov process that can realized by a stationary, reversible, continuous-time Markov process operating on each edge $e$ of $\mathcal{T}$ for a duration [corresponding to $t(e)$] (this is Theorem 4(2) of [18]; for more details on such models, see [16]).

Next, we say that a probability distribution on patterns is *separable* if it satisfies the following property:

(S) Whenever $(a_1, b_1), (a_1, b_2), \ldots, (a_m, b_m)$ are pairs of leaves whose connecting paths are pairwise edge-disjoint, then $\{E(a_i, b_i), i = 1, \ldots, m\}$ are independent events.

It is easily seen that the CFN model is separable. Moreover, any group-based model satisfies the separation condition (S) (Theorem 10 of [20], generalizing [11]); briefly, "group-based models" are defined in the same way as the CFN model, but over an arbitrary finite Abelian group, rather than the particular group $(\{0, 1\}, +_{\mathrm{mod}2})$ (for more details, see [16]).

We will call a model that satisfies conditions (C) and (S) a *conservative, separable process*. Examples of such models include the CFN model, and, more generally, the symmetric $q$-state model, for which, when a transition occurs, one of the remaining states is selected uniformly at random. For this model, we have $c = \frac{1}{q}$ in condition (C), and this model is well known in a variety of fields, including physics, broadcasting and molecular biology, where it is referred to as the "$q$-state Potts model," the "$q$-ary symmetric channel," and the "Neyman $q$-state model," respectively (and, in the special case when $q = 4$, as the "Jukes–Cantor model"); for more details, see [13]. A further example of a conservative, separable process in molecular biology is the Kimura 3ST model (for details, see [16]).

*Weakened constraints.* In Theorem 3.1 we imposed the condition $f \leq p_e$ for a fixed $f > 0$ for the transition mechanism $\mathcal{P}_2$. In fact, an inspection of the proof reveals that $0 < f = f(n)$ may depend on $n$, as far as we have $\lim_{n \to \infty} h(n) f(n) = \infty$, where $h(n)$ is any function satisfying the statement of Lemma 3.7. [The present proof of Lemma 3.7 allows $f(n) \to 0$ "very slowly," but the truth is likely just "slowly."]

The result allowing these three types of extensions is the following.

THEOREM 4.1. *Fix $0 < t_+ < \infty$, and allow $t_- = t_-(n) > 0$ to vary with $n$ if still $\lim_{n \to \infty} h(n) t_-(n) = \infty$, where $h(n)$ is any function satisfying the statement of Lemma 3.7. For every binary phylogenetic $X$-tree $\mathcal{T}_1$ with a conservative, separable process $\mathcal{P}_1$ where $t(e) \leq t_+$ in $\mathcal{P}_1$, and any $\mu$ permutation invariant measure on phylogenetic $X$-trees, the following holds for a function $\varepsilon(n) = o(1)$. The set of binary phylogenetic $X$-trees of measure $1 - o(1)$ has the property that any of them equipped with an arbitrary conservative, separable process $\mathcal{P}_2$, with $t(e) \geq t_-$ in $\mathcal{P}_2$ (assuming $\mathcal{P}_2$ has the same $H$ and $c$ as $\mathcal{P}_1$) has*

$$\text{vardist}((\mathcal{T}_1, \mathcal{P}_1), (\mathcal{T}_2, \mathcal{P}_2)) \geq 2 - \varepsilon(n). \tag{4.1}$$

PROOF. We need a straightforward modification of the proof of Theorem 3.1. Leaving out the subscript from the notation for the generic leaf pair $(a_i, b_i)$, formula (3.16) can be substituted by

$$H(3t_+) \leq H(d_{\mathcal{T}_1}(a,b) t_+) \leq \mathbb{P}[state_1(a) = state_1(b)]; \tag{4.2}$$

(3.18) can be substituted by

$$\mathbb{P}[state_2(a) = state_2(b)] \leq H(d_{\mathcal{T}_2}(a,b) t_-) \leq H(h(n) t_-) < c + \varepsilon \tag{4.3}$$

for any fixed $\varepsilon > 0$ as $n \to \infty$. For a sufficiently small $\varepsilon > 0$, we have

$$c + \varepsilon < H(3t_+) \tag{4.4}$$

[this follows from the assumptions on $H$ and $c$], and thus, inequality (4.4) substitutes for (3.20). □

**Acknowledgments.** We thank Éva Czabarka for her careful reading of earlier versions of this manuscript, and the referees for helpful comments.

VARIATIONAL DISTANCE OF TWO TREES 13

BIOMATHEMATICS RESEARCH CENTRE
DEPARTMENT OF MATHEMATICS AND STATISTICS
UNIVERSITY OF CANTERBURY
CHRISTCHURCH
NEW ZEALAND
E-MAIL: m.steel@math.canterbury.ac.nz

DEPARTMENT OF MATHEMATICS
UNIVERSITY OF SOUTH CAROLINA
COLUMBIA, SC
USA
E-MAIL: szekely@math.sc.edu